\documentclass[preprint,12pt]{elsarticle}

\usepackage{amssymb}
\usepackage{amsmath}
\usepackage{amsthm}
\usepackage{graphics,url,color,epsfig}
\usepackage{graphicx}   % for figures
\usepackage{tikz}

\usepackage[colorlinks]{hyperref}
\AtBeginDocument{%
  \hypersetup{%
    linkcolor=blue,%
    citecolor=green,%
  }%
}
\usepackage{cleveref}
\crefname{equation}{}{}

\newtheorem{theorem}{Theorem}[section]
\newtheorem{lemma}[theorem]{Lemma}

\newtheorem{remark}[theorem]{Remark}

\theoremstyle{definition}

\numberwithin{equation}{section}

\journal{***}
\begin{document}

\begin{frontmatter}

\title{Nonexistence of solutions for indefinite fractional  parabolic equations}
\author[rvt1]{Wenxiong Chen}
\ead{wchen@yu.edu}
\author[rvt2]{Leyun Wu}
\ead{leyunwu@sjtu.edu.cn}
\author[rvt3]{Pengyan Wang}
\ead{wangpy119@126.com}

\address[rvt1]{Department of Mathematical Sciences, Yeshiva University, New York, NY, 10033, USA}
\address[rvt2]{Department of Applied Mathematics, Hong Kong Polytechnic University, Hung Hom, Kowloon, Hong Kong}
\address[rvt3]{School of Mathematics and Statistics, Xinyang Normal University, Xinyang, 464000, P.R. China}

\begin{abstract}
~~~~We study  fractional  parabolic equations with  indefinite nonlinearities
$$ \frac{\partial u} {\partial t}(x,t) +(-\Delta)^s u(x,t)= x_1 u^p(x, t),\,\, (x, t) \in  \mathbb{R}^n \times \mathbb{R}, $$
where $0<s<1$ and $1<p<\infty$. We first prove that all positive bounded solutions are monotone increasing along the $x_1$ direction. Based on
this we derive a contradiction and hence obtain non-existence of solutions. These monotonicity and nonexistence results are crucial tools in a priori estimates and complete blow-up for  fractional parabolic equations in bounded domains. To this end, we introduce several new ideas and developed a systematic approach which may also be applied to investigate qualitative properties of solutions for
many other fractional parabolic problems.
\end{abstract}


\end{frontmatter}

\section{Introduction}

The fractional Laplacian  arises in a variety of
physical phenomena such as anomalous diffusion, ocean acoustic propagation, quasi-geostrophic dynamics, phase
transition, image reconstruction and so on when taking into account of the presence of long range interactions (see \cite{BG, CV, GO}). It can
also be used to model American options in mathematical finance (see \cite{A}). In recent years, considerable attentions from mathematical communities have been paid to the study of pseudo-differential equations  involving the fractional Laplacian.

 Due to the non-locality of the fractional Laplacian, many traditional approaches on local elliptic operators no longer work.
To overcome this difficulty, Caffarelli and Silvestre \cite{CS} introduced the extension method that reduced a nonlocal problem into a local one in higher dimensions. Another effective indirect approach is  the method of moving planes in integral forms \cite{CFY, DLL, LyL, LZ2}. Later, several direct methods have been developed to study fractional problems, such as the method
 of moving planes \cite{CLL, CW, DQ}, the method of moving spheres \cite{CLZr}, the sliding methods \cite{LiuC, WC2}, and others \cite{LWX, LZ1}.  Recently, for parabolic equations  involving the fractional Laplacian, Chen etc. \cite{CWNH,CWp} developed a systematical scheme to carry out the asymptotic method of moving planes to investigate qualitative properties of solutions, either on bounded or on unbounded domains.

In this paper we study  the following indefinite parabolic equations involving the fractional Laplacian:
\begin{eqnarray}\label{maineq}
\frac{\partial u} {\partial t}(x,t) +(-\Delta)^s u(x,t)= x_1 u^p(x, t),\,\, (x, t) \in  \mathbb{R}^n \times \mathbb{R},
\end{eqnarray}
where $0<s<1$ and $1<p<\infty.$

For each fixed $t\in \mathbb{R},$ the fractional Laplacian acting on $x$ is defined as
\begin{eqnarray*}
(-\Delta)^s u (x,t)&=&C_{n, s} P.V. \int_{\mathbb{R}^n} \frac{u(x, t)-u(y,t)}{|x-y|^{n+2s}}dy\\
&=&C_{n, s}  \lim_{\epsilon \rightarrow 0^{+}} \int_{\mathbb{R}^n \backslash B_\varepsilon (x)} \frac{u(x, t)-u(y,t)}{|x-y|^{n+2s}}dy,
 \end{eqnarray*}
where $P.V.$ stands for the Cauchy principal value.
It is easy to see that    for $u \in C^{1,1}_{loc}\cap   {\cal L}_{2s } ,$ $(-\Delta)^su$ is well defined, where
 $$ {\mathcal L}_{2s}=\left\{u(\cdot, t) \in L^1_{loc} (\mathbb{R}^n)\, \big| \int_{\mathbb R^n} \frac{|u(x,t)|}{1+|x|^{n+2s}}dx<+\infty\right\}.$$
 It is known that as $s \to 1,$ the fractional Laplacian $(-\Delta)^s$  goes to the regular Laplacian $-\Delta.$

We are particularly interested in solutions defined for all $t \in \mathbb{R}.$ Below we refer
them as \emph{entire solutions}.
We say that $u$ is a classical entire solution of \eqref{maineq}, if $$u(x,t)   \in \big( C^{1,1}_{loc}(\mathbb R^n)\cap {\mathcal L}_{2s} \big)\times C^1\big(\mathbb R\big) $$ and satisfies \eqref{maineq}
pointwise in  $\mathbb R^n \times \mathbb R$.
Our main goal is to prove monotonicity and Liouville theorems for solutions of \eqref{maineq}
to the effect that there are no positive bounded entire solutions, as stated in  Theorem \ref{mthm1} and Theorem \ref{mthm2} below.

It is well known that Liouville theorems  play an essential role in deriving a priori bounds for solutions to elliptic equations (\cite{BPGQ, CLL2016, FLN, GS1}) and can also be used to obtain uniqueness of solutions (\cite{C, KK, MS}).
For parabolic equations, they guarantee optimal universal estimates
of solutions of related initial and initial-boundary value problems (see \cite{PQS,Q, QS} and the
references therein).

For the indefinite elliptic  problem with the regular Laplacian:
\begin{eqnarray} \label{er1}
-\Delta u(x)=a(x_1)u^p(x) \,\mbox{ in }\, \mathbb R^n,~1<p<\infty,
\end{eqnarray}
the Liouville theorems have been extensively studied, see \cite{BCN, CL1997, DLi} and references therein.

For indefinite elliptic problems with the fractional Laplacian:
\begin{eqnarray} \label{001}
(-\Delta)^su(x)=a(x_1)f(u)\,\mbox{ in } \,\mathbb R^n,
\end{eqnarray}
Chen and Zhu \cite{CZ} obtained the nonexistence of
positive solutions to \eqref{001} for $\frac{1}{2}<s<1$  with $$a(x_1)f(u)=x_1u^p(x), ~1<p<\infty.$$ Their main approach  is to reduce the problem to a local one by means of the extension method introduced in \cite{CS}.
Subsequently, Barrios etc. \cite{BPGQ2017} generalized this result to $0<s<1$ by the method of moving planes in integral forms,
where $a$ and $f$ are nondecreasing functions satisfying some additional conditions.
At almost the same time, Chen, Li, and Zhu \cite{CLZ} proved the Liouville theorem for \eqref{001} by a direct method of moving planes under much weaker conditions than that in \cite{BPGQ2017}.

  For  indefinite  parabolic  problem with the regular Laplacian, Pol\'{a}\v{c}ik and Quittner \cite{PQ} obtained the non-existence of bounded positive solutions of the following equation:
\begin{eqnarray}\label{integerL1}
\frac{\partial u} {\partial t}(x,t) -\Delta u(x,t)= a(x_1) u^p(x,t),\,\, (x, t) \in  \mathbb{R}^n \times \mathbb{R}.
\end{eqnarray}
This Liouville theorem plays an important role in deriving a priori estimates and appears as crucial ingredients in the proof of the existence of positive equilibria which can be
efficiently used in the study of blow-up.
For instance, the nonexistence of solution for \eqref{integerL1}  can be applied to study the problem
 \begin{eqnarray}\label{reblowupproblem}
\left\{\begin{array}{ll}
 \frac{\partial u}{\partial t} (x, t) - \Delta u(x, t)= \lambda u (x,t) +x_1 u^p(x, t), & (x, t) \in  \Omega \times (0,+ \infty), \\
 u(x, t)=0, & (x, t) \in \partial \Omega \times (0,+ \infty),\\
u(x,0)=u_0(x)\geq 0,& x\in  \Omega,
\end{array} \right. \end{eqnarray}
where $\Omega$ is a smooth bounded domain in $\mathbb{R}^n,$ $\lambda \in \mathbb{R}$ and $p>1$ is subcritical.
 More specifically, the  Liouville theorem
 for \eqref{integerL1} together with known nonexistence results for the equations
\begin{eqnarray*}
\frac{\partial u} {\partial t}(x,t) -\Delta u(x,t)=  u^p(x,t),\,\, (x, t) \in  \mathbb{R}^n \times \mathbb{R},
\end{eqnarray*}
and
 \begin{eqnarray*}
 \left\{\begin{array}{ll}
\frac{\partial u} {\partial t}(x,t) -\Delta u(x,t)=  u^p(x,t),& (x, t) \in  \mathbb{R}^n_+ \times \mathbb{R},\\
u(x, t)=0,&(x, t) \in  \partial\mathbb{R}^n_+ \times \mathbb{R},
\end{array} \right.
\end{eqnarray*}
can be employed to derive suitable a priori bounds for solutions of \eqref{reblowupproblem} through blowing-up and re-scaling arguments (see \cite{BV, MZ1, PQS, Q}). And it is well known that these a priori estimates are important ingredients in obtaining the existence of
solutions of the same equations.

 As far as we know, there have been rarely any Liouville type  theorems for solutions to indefinite parabolic equations with the fractional Laplacian. The main difficulty lies in the non-locality of this operator. Many traditional approaches on local operators no longer work
 here. To circumvent these difficulties, one needs to introduce new ideas and develop new approaches.

In this paper, we establish  monotonicity and nonexistence of positive solutions for indefinite fractional parabolic equations \eqref{maineq}. Our main results are as follows.

\begin{theorem}\label{mthm1} Let $0<s<1$ and $1 < p < \infty$, and
suppose that $u(x,t)   \in \big( C^{1,1}_{loc}(\mathbb R^n)\cap {\mathcal L}_{2s} \big)\times C^1\big(\mathbb R\big) $ is a positive bounded classical solution of equation \eqref{maineq},
then  $u(x, t)$ is monotone increasing along the $x_1$ direction.
\end{theorem}

\begin{theorem}\label{mthm2}
 Let $0<s<1$ and $1 < p < \infty,$ and suppose that $u(x,t)   \in \big( C^{1,1}_{loc}(\mathbb R^n)\cap {\mathcal L}_{2s} \big)\times C^1\big(\mathbb R\big) $,  then the equation \eqref{maineq} possesses no positive bounded classical solutions.
\end{theorem}

\begin{remark}
To better illustrate the main ideas, we only consider the simple example as in equation \eqref{maineq}. The methods developed here are also applicable to deal with more general nonlinearities, for example, one can replace $x_1$ and $u^p$ in \eqref{maineq}  by $a(x_1)$ and $f(u)$ respectively  under suitable conditions and obtain the non-existence of positive solutions.

\end{remark}

To prove the above theorems, we will modify the direct method of moving planes for elliptic problems introduced in (\cite{CLL}) such that it can be applied to indefinite fractional parabolic problems. Usually, to carry out the method of moving planes, one needs to assume that the solution $u(x, t)$ tends to $0$ as $|x|\to \infty$.

For the corresponding elliptic problem without $x_1$ on the right hand side:
$$
(-\Delta)^s u(x)=u^p(x), \,\, x \in \mathbb{R}^n,
$$
without assuming $\mathop{\lim}\limits_{|x| \to \infty } u(x)=0,$ in the critical and
subcritical cases, one can exploit the Kelvin transform $v(x)=\frac{1}{|x|^{n-2s}}u(\frac{x}{|x|^2})$ to derive
\begin{eqnarray}\label{transform}
(-\Delta)^s v(x)=\frac{1}{|x|^{\gamma}}v^p(x) \mbox{ with } \mathop{\lim}\limits_{|x| \to \infty } v(x)=0.
\end{eqnarray}
Here $\gamma \geq 0$ and the coefficient $\frac{1}{|x|^\gamma}$ possesses the needed monotonicity such that one can employ the method of moving planes on the transformed equation \eqref{transform}. While in the presence of $x_1$, the coefficient of the transformed equation does not have the required monotonicity, and therefore the Kelvin transform renders useless.

However, in the process of applying this
Liouville Theorem (nonexistence of positive solutions) in the blowing up and re-scaling arguments to
establish a priori estimates, the solutions of the limiting equations are known to be only bounded. Hence when we consider equation \eqref{maineq}, it is impractical to impose the condition $\mathop{\lim}\limits_{|x| \to \infty } u(x)=0$, while it is more reasonable to assume that $u(x, t)$ is bounded. Under this weaker assumption, in order to use the method of moving planes, one needs to introduce
auxiliary functions.

To illustrate the ideas, we first set up the standard frame work in the method of moving planes.

For any given $\lambda \in \mathbb{R},$ let
$$
T_\lambda=\{ x \in \mathbb{R}^n \mid x_1 =\lambda\,\, \mbox{for}\,\, \lambda \in \mathbb{R} \}
$$
be the moving planes,
$$
\Sigma_\lambda = \{ x \in \mathbb{R}^n \mid x_1 <\lambda \}
$$
be the region to the left of the plane, and
$$
x^\lambda=(2\lambda-x_1, x_2, ..., x_n)
$$
be the reflection of $x$ about the plane $T_\lambda.$
Assume that $u(x, t)$ is a positive solution of equation \eqref{eq2}. We
compare the values of $u(x,t)$ with $u(x^\lambda, t)$ by studying the difference
$$
w_\lambda(x, t)=u(x^\lambda, t)-u(x,t).
$$

The first step in the standard method of moving planes is to show that, for $\lambda$ sufficiently negative,
$$w_\lambda (x,t) \geq 0, \;\; \forall \; (x,t) \in \Sigma_\lambda \times \mathbb{R}.$$
On the contrary, one would try to work on the negative minima of $w_\lambda$ to derive a contradiction.
However, without the decay assumption on the solutions, the minima of $w_\lambda$ may ``leak'' to infinity. To prevent this from happening, one usually constructs some auxiliary functions such as
 $$
 \bar{w}_\lambda(x, t)=\frac{w_\lambda(x, t)}{g(x)} \,\mbox{ with}\, \mathop{\lim}\limits_{|x| \to \infty } g(x)=\infty.$$
 Now $\mathop{\lim}\limits_{|x| \to \infty } \bar{w}_\lambda(x, t) =0$ and one can investigate the minima of $\bar{w}_\lambda(\cdot, t)$ for each fixed $t$.

Nonetheless,  the situation in the fractional
order parabolic equation is quite different and much more difficult than the one in the integer order parabolic equation.

When we compute the product of two functions, for the regular Laplacian,
  the convenience is:
\begin{eqnarray}\label{com}
&&\Delta_x  (\bar{w}_\lambda(x, t) g(x))\nonumber\\
&=&\Delta_x \bar{w}_\lambda(x, t) g(x)+2\nabla_x \bar{w}_\lambda(x, t)\cdot \nabla_x g(x)+ \bar{w}_\lambda (x,t) \Delta_x g(x).
\end{eqnarray}
At a minimum of $\bar{w}_\lambda(x, t),$ the middle term on the right hand side vanishes since $\nabla_x \bar{w}_\lambda(x,t)=0.$ This makes the analysis much easier.
 Pol\'{a}\v{c}ik and   Quaittner in \cite{PQ} choose $g(x)=\ln (\lambda+1-x_1)+1$ and $g(x)=2-\frac{\delta}{\delta+\lambda-x_1}$ for $x_1\leq 0$ and
 $x_1>0$ respectively.
However, the fractional counterpart of \eqref{com} is
\begin{eqnarray*}
&&(-\Delta)^s _x(\bar{w}_\lambda g)(x, t)\\
&=&(-\Delta)^s_x \bar{w}_\lambda(x, t) g(x)-2C \int_{\mathbb{R}^n} \frac{(\bar{w}_\lambda(x, t)-\bar{w}_\lambda(y, t))(g(x)-g(y))}{|x-y|^{n+2s}}dy\\
&&+\bar{w}_\lambda(x, t)(-\Delta)_x^s g(x).
\end{eqnarray*}
At a minimum of $\bar{w}_\lambda(x, t),$ the middle term on the right hand side (the integral) neither vanishes nor has a definite sign.
This makes the analysis much more difficult, and to circumvent it, as in \cite{CLZ}, the authors choose $g(x)=|x-(\lambda+1)e_1|^\sigma$ with $e_1=(1,0, \cdots, 0)$ for sufficiently small $\sigma>0$. A key estimate on $(-\Delta)^s w_\lambda$ at minima of $\bar{w}_\lambda$ plays
an important role here. However, many approaches in \cite{CLZ} for elliptic fractional equations cannot be adapted here to treat parabolic
fractional equations, and we need to introduce new ideas and develop new methods for the parabolic setting.

To prove that the solutions are increasing in $x_1$ direction, the general frame work are similar to the traditional ones:

In step 1, we show that for $\lambda \leq 0$,
\begin{equation} w_\lambda (x, t) \geq 0, \;\; \; \forall \; (x,t) \in \Sigma_\lambda \times \mathbb{R}.
\label{A11}
\end{equation}

In step 2, we move the plane toward the right as long as the above inequality holds and show
that it can be moved all the way to $x_1 = \infty$.

However, in carrying out these two steps, we employ approaches that are completely different from those in previous literature.

In step 1, to realize (\ref{A11}), we first obtain the estimate
$$
\mbox{if }\, \bar{w}_\lambda(x(t), t)<0,
\mbox{ then }\, \frac{\partial \bar {w}_\lambda}{\partial t} (x(t), {t}) \geq \frac{-C}{|x_1( t)-\lambda|^{2s}}\bar{w}_\lambda (x({t}), {t}),
$$
where
$$
\bar{w}_\lambda(x(t), t)=\inf_{x \in \Sigma_\lambda}\bar{w}_\lambda(x, t),
$$
and $x_1(t)$ is the first component of $x(t)$. Then we construct a sub-solution
$$
z(t)=-\bar{M} e^{-c_0(t-\bar{t})},
$$
and let $\bar{t} \to -\infty$ to derive a contradiction.

In step 2, let
$$
\lambda_0= \sup \{ \lambda \mid w_\mu(x, t) \geq 0, \,\, \forall (x, t) \in \Sigma_\mu \times \mathbb{R},\,\, \mu \leq \lambda\}.
$$
Suppose $\lambda_0 < \infty$. Then there exists a sequence $\lambda_k \searrow \lambda_0$,
such that $$w_{\lambda_k}(x(t_k), t_k) < 0.$$ By constructing proper
auxiliary functions, we will be able to choose a sequence of approximate minimum $(x(\bar{t}_k), \bar{t}_k)$ of $\bar{w}_{\lambda_k}$ in
$\Sigma_{\lambda_k} \times \mathbb{R}$ and derive estimates along this sequence:
$$ \frac{\partial w_{\lambda_k}}{\partial t}(x(\bar{t}_k), \bar{t}_k) \geq C_4>0.$$
While on the other hand, we show that
$$ \frac{\partial w_{\lambda_k}}{\partial t}(x(\bar{t}_k), \bar{t}_k) \to 0,$$
and hence derive a contradiction.

We believe that the new ideas introduced here can be conveniently applied to investigate qualitative properties of solutions for
other fractional parabolic problems.

In Section 2, we study the monotonicity of solutions in $x_1$-direction and establish Theorem \ref{mthm1}. In Section 3, we prove the nonexistence of positive solutions.

\section{Monotonicity}

In this section, we study the monotonicity of positive solutions for the indefinite fractional parabolic equations
\begin{eqnarray}\label{eq2}
\frac{\partial u} {\partial t}(x,t) +(-\Delta)^s u(x,t)= x_1 u^p(x, t),\,\, (x, t) \in  \mathbb{R}^n \times \mathbb{R}.
\end{eqnarray}

Employing a new version of method of moving planes, we show that the solutions are strictly monotone increasing in $x_1$-direction and
hence establish Theorem \ref{mthm1}.

For any given $\lambda \in \mathbb{R},$ let
$$
T_\lambda=\{ x \in \mathbb{R}^n \mid x_1 =\lambda\,\, \mbox{for}\,\, \lambda \in \mathbb{R} \}
$$
be the moving planes,
$$
\Sigma_\lambda:= \{ x \in \mathbb{R}^n \mid x_1 <\lambda \}
$$
be the region to the left of the plane, and
$$
x^\lambda=(2\lambda-x_1, x_2, ..., x_n)
$$
be the reflection of $x$ about the plane $T_\lambda.$
Assume that $u(x, t)$ is a positive solution of equation \eqref{eq2}. We
compare the values of $u(x,t)$ with $$u_\lambda(x, t)=u(x^\lambda, t).$$ Let
$$
w_\lambda(x, t)=u_\lambda(x, t)-u(x,t).
$$
We  deduce from the equation \eqref{eq2} that $w_\lambda$ satisfies
\begin{eqnarray}\label{eq3}
&&\frac{\partial w_\lambda} {\partial t}(x,t) +(-\Delta)^s w_\lambda(x,t)\nonumber\\
&=& x_1^\lambda u_\lambda^p(x, t)-x_1u^p(x, t)\nonumber\\
&=& (x_1^\lambda-x_1)u_\lambda^p(x, t)+x_1(u_\lambda^p(x, t)-u^p(x,t))\nonumber\\
&\geq & x_1(u_\lambda^p(x, t)-u^p(x,t))\nonumber\\
&=& x_1 p \xi_\lambda^{p-1}(x, t) w_\lambda(x, t)\nonumber\\
&:=& c_\lambda(x,t) w_\lambda(x,t),
\end{eqnarray}
where $\xi_\lambda(x, t)$ lies between $u_\lambda(x, t)$ and $u(x,t),$
and the sign of $c_\lambda(x,t)$   is consistent with $x_1$ and hence with $\lambda$ in $\Sigma_\lambda \times \mathbb{R}.$
Therefore it is reasonable to handle the  two cases:
$$ c_\lambda(x, t)\leq 0 \,\, \; \mbox{ and } \;\;   c_\lambda(x, t)> 0$$
separately. We divide the proof of Theorem \ref{mthm1} into two steps.

In the first step, we show that for any $\lambda \leq 0$,
\begin{equation}
w_\lambda(x, t) \geq 0 \mbox{ in } \Sigma_\lambda \times \mathbb{R}.
\label{A3}
\end{equation}
In this case, the crucial fact that $c_\lambda (x,t) \leq 0$ enables us to construct a sub-solution $z(t)$ in a certain cylinder
to derive (\ref{A3}), which provides a starting point to move the plane.

In Step 2, we move the plane $T_\lambda$ toward the right as long as inequality (\ref{A3}) holds to its limiting position
$T_{\lambda_0}$. We show that $\lambda_0$ must be $\infty$.

\medskip

{\bf{Proof of Theorem \ref{mthm1}.} }

Recall the problem we considered (see \eqref{eq3}):
\begin{eqnarray}\label{wlambda}
\left\{\begin{array}{ll}
 \frac{\partial w_\lambda}{\partial t} (x, t)+( - \Delta )^s w_\lambda(x, t)\geq c_\lambda (x,t) w_\lambda (x,t), & (x, t) \in  \Sigma_\lambda \times \mathbb{R}, \\
w_\lambda(x^\lambda,t) =-w_\lambda(x, t),& (x,t )\in  \Sigma_\lambda \times \mathbb{R}.
\end{array} \right. \end{eqnarray}

\emph{Step 1.}
In this step, we show that for $\lambda \leq 0,$ we have
\begin{eqnarray}\label{eq-1}
w_\lambda(x, t) \geq 0, \,\, (x, t) \in  \Sigma_\lambda \times \mathbb{R},
\end{eqnarray}
which will be proved by adopting the idea of Theorem 5 in \cite{CW1}.

If $\lambda \leq 0, $ by the definition of $c_\lambda(x, t)$ in \eqref{eq3}, one has
$$c_\lambda(x, t) \leq 0 \mbox{ in }  \Sigma_\lambda \times \mathbb{R}.$$

Since we do not impose any decay condition on $w_\lambda$, to prevent its minimum from leaking to infinity, we introduce an auxiliary function
\begin{equation}
g(x)=|x-(\lambda+1)e_1|^\sigma,
\label{A2}
\end{equation}
where $e_1=(1, 0, ...,0)$ and $\sigma$ is a small positive number to be chosen as in the proof of Theorem 1 in \cite{CLZ}.
We consider
$$
 \bar{w}_\lambda(x, t)=\frac{w_\lambda(x, t)}{g(x)}.$$
Obviously, $\bar{w}_\lambda(x, t)$ and $w_\lambda(x, t)$ have the same sign, however
  \begin{equation}
  \lim_{|x|\to +\infty} \bar{w}_\lambda(x, t) =0.
  \label{A1}
  \end{equation}
  Therefore, we will consider the function $\bar{w}_\lambda(x, t)$  in subsequent processes.

For any  fixed
$t\in \mathbb{R},$ (\ref{A1}) implies that there exists $x(t)$, such that
$$
\bar{w}_\lambda(x(t), t)=\inf_{x \in \Sigma_\lambda}\bar{w}_\lambda(x, t).
$$

First, we conclude that for any  fixed
$t\in \mathbb{R}$ and $\lambda \leq 0,$
\begin{eqnarray}\label{keyineq}
\mbox{if }\, \bar{w}_\lambda(x(t), t)<0,
\mbox{ then }\, \frac{\partial \bar {w}_\lambda}{\partial t} (x(t), {t}) \geq \frac{-C}{|x_1( t)-\lambda|^{2s}}\bar{w}_\lambda (x({t}), {t}).
\end{eqnarray}

In fact, by a similar calculation as (22) in \cite{CLZ}, we derive that $$\mbox{ if }  w_\lambda(x(t), t)<0, \mbox{ then }
(-\Delta)^s w_\lambda (x(t), t) \leq \frac{C}{|x_1(t)-\lambda|^{2s}}w_\lambda (x(t), t).
$$
Combining this with \eqref{wlambda}, one has
$$
\frac{\partial w_\lambda}{\partial t} (x(t), t) \geq \frac{-C}{|x_1(t)-\lambda|^{2s}}w_\lambda (x(t), t).
$$
Then by the definition of  $\bar{w}_\lambda (x, t),$ we derive \eqref{keyineq}.
\medskip

For any  fixed $t\in \mathbb{R},$ define
$$
m(t):=\bar{w}_\lambda(x(t), t)=\inf_{x \in \Sigma_\lambda}\bar{w}_\lambda(x, t).
$$
In order to prove \eqref{eq-1}, it suffices to show that
\begin{eqnarray}\label{mtgeq}
 m(t)\geq 0, \, \forall~ t \in \mathbb{R}.
\end{eqnarray}

To proceed with the proof, we need the following lemma.
\medskip

\begin{lemma}\label{conslem}
Given any $m_0>0$, there exists a positive constant $c_0$ such that  if $m(t)\leq -m_0,$ then
\begin{eqnarray}\label{c0}
\frac{C}{|x_1(t)-\lambda|^{2s}} >c_0>0,
\end{eqnarray}
where $x(t)=(x_1(t), ... , x_n(t))$ is a minimum point of $\bar{w}_\lambda (x, t)$ in $\Sigma_\lambda$ for each fixed $t$.
\end{lemma}
{\bf{Proof.}} If the conclusion of Lemma \ref{conslem} is not valid, then there exist an $\bar{m}_0 >0$ and  a sequence of $\{t_k\},\, k=1,2, ...$ such that
\begin{eqnarray}\label{ck}
m(t_k) \leq - \bar{m}_0,
\end{eqnarray}
and
$$
\frac{C}{|x_1(t_k)-\lambda|^{2s}} \to 0, \,\, k\to +\infty,
$$
therefore,
$$
|x(t_k) | \to +\infty,\,\, k\to +\infty,
$$
and it follows that
$$
m(t_k) =\bar {w}_\lambda (x(t_k), t_k) \to 0,\,\, k\to +\infty.
$$
This contradicts \eqref{ck} and thus completes the proof of Lemma \ref{conslem}.
\medskip

Now we continue the proof of \eqref{mtgeq}.

If \eqref{mtgeq} is false, then  there exits  a $t_0 \in \mathbb{R}$ such that
\begin{eqnarray}\label{keystep1}
-m_0:=m(t_0)=\bar{w}_\lambda(x(t_0), t_0)<0.
\end{eqnarray}

For any $\bar{t}<t_0,$ we construct a subsolution as
$$
z(t)=-\bar{M} e^{-c_0(t-\bar{t})},
$$
where $c_0$ is as defined in \eqref{c0} and
$$
-\bar{M}= \inf_{\Sigma_\lambda \times \mathbb{R}} \bar{w}_\lambda(x, t).
$$

We will prove that
\begin{eqnarray}\label{subsolu}
\bar{w}_\lambda(x, t) \geq z(t) \, \mbox{ in the parabolic cylinder } \overline{\Sigma_\lambda} \times [ \bar{t}, t_0].
\end{eqnarray}
Consider the function
$$
v(x, t)= \bar{w}_\lambda(x, t) - z(t), \,\, (x, t) \in \overline{\Sigma_\lambda} \times [ \bar{t}, t_0].
$$
By our construction of $z(t)$, we  have on the bottom of the cylinder,
$$
v(x, t) =\bar{w}_\lambda(x, t)-z(t) =\bar{w}_\lambda(x, t)-(-\bar{M}) \geq 0,\,\, (x, t) \in \Sigma_\lambda \times \{\bar{t}\};
$$
while on the side of the cylinder, we also have
$$
v(x, t) =\bar{w}_\lambda(x, t)-z(t) =-z(t) \geq 0,\,\, (x, t) \in T_\lambda \times [ \bar{t}, t_0].
$$

See Figure 1 below.
\bigskip

\begin{center}
\begin{tikzpicture}[node distance = 0.3cm]
\draw [->, semithick] (-3.55,0) -- (4,0) node[right] {$x_1$};

 \fill [blue!50] (1,0,0)--( 1,2 ,0)--( 1,2,-2)
                                         --(1,0,-2)--cycle;
 \path (1,0) [very thick,fill=black]  circle(1.5pt) node at (1,-0.3) {$\bar t$};

 \path (-1,3.4) node at (-0.6,3.4) {$ \Sigma_\lambda\times [\bar t,t_0]$};
 \fill [orange] (-3.55,0,0)--( 0.26,0 ,0)--( -0.53,0,-2) --(-3.54,0,-2)--cycle;
 \fill [orange] (-3.55,0,0.03)--( 0.26,0 ,0.03)--( 1,0,2) --(-3.6,0,2)--cycle;
 \fill [blue!25] (-3.55,2,0)--( 1,2 ,0)--( 1,2,-2) --(-3.54,2,-2)--cycle;
 \fill [blue!25] (-3.55,2,-0.1)--( 1,2 ,-0.1)--( 1,2,2) --(-3.6,2,2)--cycle;
 \fill [blue!50](1,0,0)--( 1,2 ,0)--( 1,2,2)
                                         --(1,0,2)--cycle;
 \path (1,2,0) [very thick,fill=black]  circle(1pt) node at (1.2,1.8,0) {$ t_0$};
 \draw [->, semithick] (1,0) -- (1,3.3,0) node[right] {$t$};
\draw [->, semithick] (1, 0,2) -- (1,0,-2) node[right] {$T_\lambda$};
 \path (-2.6 ,0,-1) node [black]{$\Sigma_\lambda$};
\node [below=0.8cm, align=flush center,text width=8cm] at (0,-0.4)
        {$Fig.$1.    \fontsize{10}{10}\selectfont  {The cylinder} };
\end{tikzpicture}
\end{center}

It follows that if \eqref{subsolu} is false, then there exists a point $(x(\tilde{t}), \tilde{t}) \in  \Sigma_\lambda \times (\bar{t}, t_0]$ such that
\begin{eqnarray}\label{minv}
v (x(\tilde{t}), \tilde{t}) =\inf_{\overline{\Sigma_\lambda} \times [\bar{t}, t_0]} v(x, t)<0.
\end{eqnarray}
Obviously,
\begin{eqnarray}\label{minvt}
\frac{\partial v}{\partial t} (x(\tilde{t}), \tilde{t}) \leq 0.
\end{eqnarray}

On one hand, from the definition of $v(x,t),$ we have
$$
\bar{w}_\lambda(x(\tilde{t}), \tilde{t})=\inf_{{\Sigma_\lambda}} \bar{w}_\lambda(x, \tilde{t})< z (\tilde t)<0.
$$
Therefore, by \eqref{keyineq}, we have
\begin{eqnarray}\label{cc}
\frac{\partial \bar {w}_\lambda}{\partial t} (x(\tilde{t}), \tilde{t}) \geq \frac{-C}{|x_1(\tilde t)-\lambda|^{2s}}\bar{w}_\lambda (x(\tilde{t}), \tilde{t}).
\end{eqnarray}

On the other hand, we obtain from \eqref{minv} that
$$
v (x(\tilde{t}), \tilde{t}) \leq v (x(t_0), t_0),
$$
i.e.,
$$
\bar{w}_\lambda(x(\tilde{t}), \tilde{t})-\bar{w}_\lambda(x(t_0), t_0)\leq z(\tilde{t})-z(t_0)\leq 0
$$
due to the monotonicity of $z(t).$ Therefore,
\begin{eqnarray}\label{lem3condition}
m(\tilde{t})= \bar{w}_\lambda(x(\tilde{t}), \tilde{t}) \leq \bar{w}_\lambda(x(t_0), t_0)=m(t_0)=-m_0.
\end{eqnarray}
Using Lemma \ref{conslem}, we derive from \eqref{lem3condition} that
$$
\frac{C}{|x_1(\tilde{t})-\lambda|^{2s}} >c_0>0.
$$
Combining this with \eqref{cc}, we obtain
\begin{eqnarray}\label{posic}
\frac{\partial \bar {w}_\lambda}{\partial t} (x(\tilde{t}), \tilde{t}) \geq -c_0\bar{w}_\lambda (x(\tilde{t}), \tilde{t}).
\end{eqnarray}
Then by \eqref{minvt}, we derive
$$
-c_0 z(\tilde t)= \frac {\partial z} {\partial t}(\tilde t) \geq \frac {\partial \bar{w}_\lambda} {\partial t} (x(\tilde t), \tilde t) \geq -c_0 \bar{w}_\lambda (x(\tilde t), \tilde t),
$$
it follows that
$$
v(x(\tilde{t}), \tilde{t})=\bar{w}_\lambda(x(\tilde t), \tilde t)-z(\tilde t) \geq 0,
$$
which  contradicts
$$
v (x(\tilde{t}), \tilde{t})<0.
$$
Therefore, we conclude that \eqref{subsolu} holds. It follows that
$$
\bar{w}_\lambda (x, t)\geq z(t),\,\, (x, t) \in \overline{\Sigma_\lambda} \times [\bar{t},  t_0].
$$
Since the above is true for any $\bar{t}$, we can let   $\bar{t} \to -\infty$ to conclude that $ z(t) \to 0$, and hence
$$
\bar{w}_\lambda (x, t) \geq 0,\,\, (x, t) \in \overline{\Sigma_\lambda} \times (-\infty,  t_0].
$$
This contradicts \eqref{keystep1}. Therefore \eqref{mtgeq} must be valid and so does \eqref{eq-1}.
\medskip

\begin{remark}\label{remarkstep1}
By a similar proof as in Step 1, we can derive that for any $\lambda>0$,  if $w_\lambda(x, t)$ is a solution of \eqref{wlambda} with
$$
\bar{w}_\lambda(x(t), t)=\inf_{x \in \Sigma_\lambda}\bar{w}_\lambda(x, t) <0,
$$
then $x_1(t)>0.$ This conclusion will be used in Step 2.
\end{remark}

\emph{Step 2.}
In this step, we  move the plane to the right as long as
 \begin{eqnarray}\label{eq2-1}
w_\lambda(x, t) \geq 0, \,\, (x, t) \in  \Sigma_\lambda \times \mathbb{R}
\end{eqnarray}
holds to its limiting position.

Define
$$
\lambda_0= \sup \{ \lambda \mid w_\mu(x, t) \geq 0, \,\, \forall (x, t) \in \Sigma_\mu \times \mathbb{R},\,\, \mu \leq \lambda\}.
$$
We show that
\begin{eqnarray}\label{eq2-2}
\lambda_0=+\infty.
\end{eqnarray}

Suppose in the contrary, $0<\lambda_0<+\infty,$ we will derive a contradiction.

By the definition of $\lambda_0$, there exists a sequence $\lambda_k \searrow \lambda_0,$ such that
$$\inf_{\Sigma_{\lambda_k} \times \mathbb{R}} w_{\lambda_k}(x, t)<0.$$
Denote
$$
 \bar{w}_{\lambda_k}(x, t)=\frac{w_{\lambda_k}(x, t)}{g(x)}$$
 with $g(x)$ defined in (\ref{A2}).
 Then obviously,
\begin{eqnarray}\label{eq2-3}
-m_k:= \inf_{\Sigma_{\lambda_k} \times \mathbb{R}} \bar{w}_{\lambda_k}(x, t)<0.
\end{eqnarray}
Since $\mathbb{R}$ is unbounded, the minimum value of $\bar{w}_{\lambda_k}$ may not be attained for some finite value $t$. In order to obtain
more information on $\frac{\partial \bar{w}_{\lambda_k}}{\partial t}$,  we choose a sequence $t_k$, and corresponding $x(t_k)$ and $\varepsilon_k \searrow 0$ such that
\begin{eqnarray}\label{eq2-4}
\bar{w}_{\lambda_k}(x(t_k), t_k)=\inf_{\Sigma_{\lambda_k}}\bar{w}_{\lambda_k}(\cdot, t_k)=-m_k+\varepsilon_k m_k.
\end{eqnarray}
We introduce an auxiliary function
\begin{eqnarray*}
\tilde{w}_{\lambda_k}(x, t)=\bar{w}_{\lambda_k}(x, t)-\varepsilon_k m_k\eta_k(t),
\end{eqnarray*}
where $\eta_k(t)=\eta(t-t_k),$ $\eta(t)\in C_0^\infty(\mathbb R)$, $|\eta'(t)|\leq 1$  and
\begin{eqnarray*}
 \eta(t)=\left\{\begin{array}{ll}
1,& \quad  |t| \leq \frac{1}{2},\\
 0,&\quad  |t| \geq 2 .
\end{array}
\right.
\end{eqnarray*}
Now the  minimum value  of $\tilde{w}_{\lambda_k}(x, t)$ in $\Sigma_{\lambda_k} \times (t_k-2, t_k+2)$ is less or equal to
the value in its complement due to the truncation of $\eta_k(t)$. Therefore, the minimum of $\tilde{w}_{\lambda_k}(x, t)$ in  $\Sigma_{\lambda_k} \times \mathbb{R}^n$ is attained in $\Sigma_{\lambda_k} \times (t_k-2, t_k+2),$ along which we will be able to derive a contradiction. More precisely,
by the definition of $\tilde{w}_{\lambda_k}(x, t),$ one has
$$
\tilde{w}_{\lambda_k}(x(t_k), t_k)=-m_k,
$$
however, for $|t-t_k|\geq 2,$
$$
\tilde{w}_{\lambda_k}(x,t)=\bar{w}_{\lambda_k}(x,t) \geq -m_k.
$$
Therefore,
$\tilde{w}_{\lambda_k}(x,t)$ attains its minima in $\Sigma_{\lambda_k} \times (t_k-2, t_k+2),$ say at $(x(\bar{t}_k), \bar{t}_k),$ i.e.,
$$
\tilde{w}_{\lambda_k}(x(\bar{t}_k), \bar{t}_k)=\inf_{\Sigma_{\lambda_k} \times \mathbb{R}}\tilde{w}_{\lambda_k}(x, t)<0.
$$
Then
$$
\frac{\partial \tilde{w}_{\lambda_k}}{\partial t}(x(\bar{t}_k), \bar{t}_k)=0,
$$
and it follows that
\begin{eqnarray}\label{eq2-5}
\left|  \frac{\partial \bar{w}_{\lambda_k}}{\partial t}(x(\bar{t}_k), \bar{t}_k) \right|  = \left| \varepsilon_k m_k \frac{\partial \eta_k}{\partial t}\right| \leq \varepsilon_k m_k.
\end{eqnarray}

From the definition of $m_k$ in \eqref{eq2-3} and
$$
\tilde{w}_{\lambda_k}(x(\bar{t}_k), \bar{t}_k) \leq \tilde{w}_{\lambda_k}(x(t_k), t_k),
$$
we obtain
\begin{eqnarray}\label{eq2-6}
-m_k\leq \bar{w}_{\lambda_k}(x(\bar{t}_k), \bar{t}_k) \leq \bar{w}_{\lambda_k}(x(t_k), t_k)=-m_k+\varepsilon_k m_k.
\end{eqnarray}
By the definition of $\tilde{w}_{\lambda_k}(x, t),$ we know that
$$
\bar{w}_{\lambda_k}(x(\bar{t}_k), \bar{t}_k)=\inf_{x \in \Sigma_{\lambda_k}}\bar{w}_{\lambda_k}(x, \bar{t}_k).
$$

From Remark \ref{remarkstep1}, we derive that $x_1(\bar{t}_k) >0.$ Hence we may assume $0<x_1(\bar{t}_k)<\lambda_0+1$. Then by a similar calculation as (22) in \cite{CLZ}, we obtain
\begin{eqnarray}\label{cank}
(-\Delta)^s w_{\lambda_k} (x(\bar{t}_k), \bar{t}_k) \leq \frac{C}{|x_1(\bar{t}_k)-\lambda_k|^{2s}}w_{\lambda_k} (x(\bar{t}_k), \bar{t}_k).
\end{eqnarray}

Notice that there exists a positive constant $C_1$ such that
$$
x_1(\bar{t}_k) p \xi_{\lambda_k}^{p-1}(x(\bar{t}_k), \bar{t}_k)\leq C_1.
$$
By \eqref{wlambda}, \eqref{cank} and  $w_{\lambda_k} (x(\bar{t}_k), \bar{t}_k)<0$, we have
\begin{eqnarray}\label{initial}
&&\frac{\partial w_{\lambda_k}}{\partial t}(x(\bar{t}_k), \bar{t}_k)+\frac{C}{|x_1(\bar{t}_k)-\lambda_k|^{2s}}w_{\lambda_k} (x(\bar{t}_k), \bar{t}_k)\nonumber\\
&\geq&
x_1(\bar{t}_k) p \xi_{\lambda_k}^{p-1}(x(\bar{t}_k), \bar{t}_k) w_{\lambda_k}(x(\bar{t}_k), \bar{t}_k)\nonumber\\
&\geq  & C_1 w_{\lambda_k} (x(\bar{t}_k), \bar{t}_k).
\end{eqnarray}
Again dividing $g(x(\bar{t}_k), \bar{t}_k)$ on both sides of the above inequality, we obtain
\begin{eqnarray}\label{eq2-7}
\frac{\partial \bar{w}_{\lambda_k}}{\partial t}(x(\bar{t}_k), \bar{t}_k)+\frac{C}{|x_1(\bar{t}_k)-\lambda_k|^{2s}}\bar{w}_{\lambda_k} (x(\bar{t}_k), \bar{t}_k)\geq C_1 \bar{w}_{\lambda_k} (x(\bar{t}_k), \bar{t}_k).
\end{eqnarray}
Combine \eqref{eq2-5}, \eqref{eq2-6} and \eqref{eq2-7}, then divide $-m_k$ on both sides, since $\varepsilon_k$ is small, one can show that
\begin{eqnarray}\label{eq2-8}
\frac{C}{|x_1(\bar{t}_k)-\lambda_k|^{2s}} \leq \frac{C_1}{2}.
\end{eqnarray}
It follows that if $k$ is large enough, we have
$$
|x_1(\bar{t}_k)-\lambda_k| \geq C_2>0
$$
and
\begin{eqnarray}\label{eq2-9}
|x_1(\bar{t}_k)-\lambda_0| \geq \frac{C_2}{2}>0.
\end{eqnarray}

More accurately, from the initial inequality \eqref{initial}, we are able to modify \eqref{eq2-8} as
\begin{eqnarray*}
\frac{C}{|x_1(\bar{t}_k)-\lambda_k|^{2s}}
&\leq& \frac{p}{2}  x_1(\bar{t}_k) \xi^{p-1}_{\lambda_k}(x(\bar{t}_k), \bar{t}_k) \\
&\leq& \frac{p}{2}  x_1(\bar{t}_k) u^{p-1}(x(\bar{t}_k), \bar{t}_k)\\
&\leq& \frac{p}{2}  \lambda_k u^{p-1}(x(\bar{t}_k), \bar{t}_k).
\end{eqnarray*}
Notice that $|x_1(\bar{t}_k)-\lambda_k|\leq \lambda_k \leq \lambda_0+1$ for sufficiently large $k$, therefore,  there exists a positive constant $C_3$ such that
\begin{eqnarray}\label{eq2-10}
u(x (\bar{t}_k), \bar{t}_k) \geq C_3>0.
\end{eqnarray}

By  \eqref{eq3}, \eqref{eq2-9} and \eqref{eq2-10}, for $k$ sufficiently large, we have
\begin{eqnarray}\label{eq2-11}
&&\frac{\partial w_{\lambda_k}}{\partial t}(x(\bar{t}_k), \bar{t}_k)+\frac{C}{|x_1(\bar{t}_k)-\lambda_k|^{2s}}w_{\lambda_k} (x(\bar{t}_k), \bar{t}_k)\nonumber\\
&\geq & (x_1^\lambda(\bar{t}_k)-x_1(\bar{t}_k)) u^p_{\lambda_k} (x(\bar{t}_k), \bar{t}_k)+x_1(\bar{t}_k) p u^{p-1}(x(\bar{t}_k), \bar{t}_k) w_{\lambda_k}(x(\bar{t}_k), \bar{t}_k)\nonumber\\
&\geq& C_4>0,
\end{eqnarray}
where we have used the fact that $u_{\lambda_k}(x, t)$ is uniformly H\"{o}lder continuous and
$$
u_{\lambda_k}(x, t) \rightrightarrows u_{\lambda_0}(x, t) \,\, \mbox{and}\,\, w_{\lambda_k}(x, t) \rightrightarrows w_{\lambda_0}(x, t) \geq 0.
$$
Since $w_{\lambda_k} (x(\bar{t}_k), \bar{t}_k)<0,$ we derive from \eqref{eq2-11} that
\begin{eqnarray}\label{eq2-12}
\frac{\partial w_{\lambda_k}}{\partial t}(x(\bar{t}_k), \bar{t}_k) \geq C_4>0.
\end{eqnarray}

Next we translate $w_{\lambda_k}$, let
$$
\hat w_{\lambda_k} (x ,t)= w_{\lambda_k} (x+x(\bar{t}_k), t+\bar{t}_k),
$$
we obtain from \eqref{eq2-12} that
\begin{eqnarray}\label{eq2-13}
\frac{\partial \hat w_{\lambda_k}}{\partial t}(0, 0) \geq C_4>0.
\end{eqnarray}
By Theorem 1.3 and Theorem 1.1 in \cite{FR}, we derive that
$$
\| \hat w_{\lambda_k}  \|_{C_{t, x}^{1+\varepsilon, 2s(1+\varepsilon)}} \leq C, \,\, \forall~ (x, t) \in \Omega \times (-T, T)\subset \mathbb{R}^n \times \mathbb{R},
$$
it follows that there exists a subsequence of $(x(\bar{t}_k), \bar{t}_k)$ (we still denote it as $(x(\bar{t}_k), \bar{t}_k)$) such that as $k \to +\infty,$
$$
\hat w_{\lambda_k} (x ,t)\to \hat w_{\lambda_0} (x ,t) \,\, \mbox{and}\,\, \frac{\partial \hat w_{\lambda_k}}{\partial t} (x ,t)\to \frac{\partial \hat w_{\lambda_0}}{\partial t} (x ,t).
$$

Since
$$
0< x_1(\bar{t}_k)\leq \lambda_k,
$$
and
$$
\lambda_k \to \lambda_0 \, \,\mbox{ as }\,\, k\to +\infty.
$$
Therefore, there exists a subsequence of $x_1(\bar{t}_k)$ (still denoted by  $x_1(\bar{t}_k)$) and $0\leq x_1^0 \leq \lambda_0$ such that
$$
 x_1(\bar{t}_k) \to x_1^0.
$$

Now consider the  function $\hat w_{\lambda_0}(x, t).$  Obviously,  we have
$$
\hat w_{\lambda_0}(x, t) \geq 0, \,\, (x, t) \in \Sigma_{\lambda_0-x_1^0} \times \mathbb{R}.
$$
Since
$$
w_{\lambda_k} (x(\bar{t}_k), \bar{t}_k)<0,
$$
we derive
$$
\hat w_{\lambda_0}(0, 0)= 0=\inf_{\Sigma_{\lambda_0-x_1^0} \times \mathbb{R}}\hat w_{\lambda_0}(x, t).
$$
Then
$$
\frac{\partial \hat w_{\lambda_0}}{\partial t}(0, 0)= 0.
$$
This contradicts  \eqref{eq2-13}, therefore, we must have
$
\lambda_0=+\infty.
$
\medskip

Therefore, $u(x, t)$ is monotone increasing along the $x_1$-direction.

This completes the proof of Theorem \ref{mthm1}.
\medskip

\section{Nonexistence}

In the proof of Theorem \ref{mthm1}, we have shown that positive solutions of \eqref{eq2}
is monotone increasing along the $x_1$ direction. Based on this, we will derive a contradiction to obtain the non-existence of solutions to \eqref{eq2} and hence prove Theorem \ref{mthm2},
by virtue of the first eigenfunction, mollification technique, and an integration by parts inequality for the nonlocal fractional Laplacian.
\medskip

{\bf{Proof of Theorem \ref{mthm2}.} }

We use a contradiction argument. Assume that there exists a  positive bounded solution $u$ of \eqref{eq2}, we will derive a contradiction.

Let $\lambda_1$ be the first eigenvalue of the problem
\begin{eqnarray*}
\left\{\begin{array}{ll}
 ( - \Delta )^s \phi (x) =\lambda_1  \phi (x),& x \in B_1(a+2, 0'),\\
\phi (x)=0,& x \in B_1^c(a+2, 0'),
\end{array} \right.
\end{eqnarray*}
where $0\leq a\in \mathbb{R}$ will be chosen sufficiently large.

In order to do integration by parts, we mollify it to be
$$\phi_1(x)=\rho \ast \phi(x) \in C_0^\infty (\mathbb{R}^n).$$ Then
\begin{eqnarray}\label{Appen}
 ( - \Delta )^s \phi_1 (x) \leq \lambda_1  \phi_1 (x), \,\, x \in \mathbb{R}^n,
\end{eqnarray}
where $\ast$ denotes convolution, $\rho(x) \in C_0^\infty (B_1(a+2, 0'))$ is a mollifier satisfying
$\int_{\mathbb{R}^n}\rho(x)dx =1.$

Eq. \eqref{Appen} actually is an integration by parts inequality, whose proof adopts  the idea of   \cite[Lemma 2.3]{LWX}, and will be presented in the Appendix (see Lemma \ref{mollification}).

We may assume that
$$
\int_{\mathbb{R}^n} \phi_1 (x) dx =1.
$$

The support of $\phi_1$ is contained in $B_2(a+2, 0')$ due to the mollification. Set
$$
\psi_a(t):= \int_{\mathbb{R}^n} u(x, t)\phi_1 (x)dx = \int_{B_2(a+2, 0')} u(x, t)\phi_1 (x)dx.
$$

By Jensen inequality, Remark \ref{mrem-inter} and \eqref{Appen}, we conclude that,
\begin{eqnarray}\label{con}
\frac{d }{dt}\psi_a(t)&=& -\int_{\mathbb{R}^n}(-\Delta)^su(x,t)\phi_1 (x)dx + \int_{\mathbb{R}^n} x_1 u^p(x,t)\phi_1 (x)dx\nonumber\\
&=&-\int_{\mathbb{R}^n}u(x,t)(-\Delta)^s \phi_1 (x)dx + \int_{\mathbb{R}^n} x_1 u^p(x,t)\phi_1 (x)dx\nonumber\\
&\geq&- \lambda\int_{\mathbb{R}^n}u(x,t)\phi_1 (x)dx + a \int_{\mathbb{R}^n}  u^p(x,t)\phi_1 (x)dx\nonumber\\
&\geq & - \lambda \psi_a(t) +  a \left( \int_{\mathbb{R}^n}  u(x,t)\phi_1 (x)dx\right)^p\nonumber\\
&= & - \lambda \psi_a(t) +  a \psi_a^{p-1}(t) \psi_a(t)\nonumber\\
&=&(a \psi_a^{p-1}(t)-\lambda)\psi_a(t).
\end{eqnarray}
 Since $u(x, t)$ is increasing in $x_1$ by Theorem \ref{mthm1}, then for any fixed $t \in \mathbb{R},$ $\psi_a(t)$ is monotone increasing with respect to $a.$
 Therefore,
\begin{eqnarray}\label{psi0}
\psi_a(0)\geq 2c_0:=\psi_0(0).
\end{eqnarray}
If $t \geq 0$ is such that $\psi_a(t) \geq c_0,$, we can choose $a$ large enough so that
$$a \psi_a^{p-1}(t)-\lambda \geq 1, $$
then it follows from \eqref{con} that
$$
\frac{d }{dt}\psi_a(t) \geq  \psi_a(t) .
$$
Thus we deduce from \eqref{psi0}  that
$$
\psi_a(t) \geq 2c_0 e^t.
$$
Now we verify  the condition that
\begin{eqnarray}\label{condition}
\psi_a(t) \geq c_0, \,\, \forall ~t \geq 0
\end{eqnarray}
by a contradiction argument.
Suppose \eqref{condition} is  false, then there exists $t_0 >0$ such that
$$
\psi_a(t_0)=c_0 \,\, \mbox{ and }\,\, \psi_a(t)>c_0 \,\, \mbox{ in }\,\, [0, t_0).
$$
It follows that $\psi_a(t) \geq 2c_0 e^t \geq 2c_0$ in $[0, t_0),$  this contradicts $\psi_a(t_0)=c_0$ and we derive \eqref{condition}.

Therefore, $\psi_a(t)$ is monotone increasing with respect to $t,$ and
 $$
 \psi_a(t) \geq 2c_0 e^t,\,\, \forall~ t \geq 0.
 $$
Consequently,
$$
\psi_a(t) \to +\infty, \,\, \mbox{as} \,\, t\to +\infty,
$$
which contradicts the boundedness of $u(x, t).$

This completes the proof of Theorem \ref{mthm2}.

\section{Appendix}

\begin{lemma}\label{mollification}
Denote by $\phi$ the first eigenfunction associated with $(-\Delta)^s$ in $B_1(0):$
\begin{eqnarray*}
\left\{\begin{array}{ll}
 ( - \Delta )^s \phi (x) =\lambda_1  \phi (x),& x \in B_1(0),\\
\phi (x)=0,& x \in B^c_1(0).
\end{array} \right.
\end{eqnarray*}
Let $\rho(x) \in C_0^\infty (B_1(0))$  and satisfy $\int_{B_1(0)}\rho(x)dx =1.$
Then we have
\begin{eqnarray}\label{eq2-mthm}
\int_{\mathbb{R}^n} (-\Delta)_z^s \phi(z) \rho(x-z)dz =\int_{\mathbb{R}^n} \phi(z)  (-\Delta)_z^s \rho (x-z)dz,
\end{eqnarray}
and
\begin{eqnarray}\label{inter}
( - \Delta )^s \phi_1 (x) \leq \lambda_1  \phi_1 (x),\,\, x \in \mathbb{R}^n,
\end{eqnarray}
where the mollification $\phi_1 (x)=(\phi \ast \rho) (x),$ and $\ast$ denotes the convolution.
\end{lemma}

{\bf{Proof.} }

\emph{Step 1.}
In this step, we estimate  $(-\Delta)^s\phi (x)$ in $ B_1^c(0),$ i.e., we show that
\begin{eqnarray}\label{eq2-16}
\left| (-\Delta)^s \phi(x) \right| \sim \frac{1}{dist^s(x, \partial B_1)},\,\, x\in B_1^c(0).
\end{eqnarray}

From \cite{QX,SV1,SV2}, one knows that the eigenvalue function $\phi (x)\in C^{1,1}_{loc}\cap {\mathcal  L}_{2s}$.
If $x \in B_1^c(0), $ by the global $s$-H\"{o}lder continuity, we have
\begin{eqnarray}\label{eq2-14}
 \left | (-\Delta)^s \phi(x)\right|
&=&C_{n, s}\left | \int_{B_1(0)}\frac{\phi(y)}{|x-y|^{n+2s}}dy\right|\nonumber\\
&\leq &C \int_{B_1(0)}\frac{(dist (y, \partial B))^s}{|x-y|^{n+2s}}dy\nonumber\\
&\leq &C \int_{B_1(0)}\frac{1}{|x-y|^{n+s}}dy.
\end{eqnarray}

We estimate the integral in the last line of \eqref{eq2-14}.

Denote $x=(x_1, 0'), \,\, x_1<0,$ and
$$
D=\{y \mid 0 < y_1< 2,\,\, |y'|<1 \};
$$
\begin{center}
\begin{tikzpicture}[node distance = 0.3cm]
\draw [->, semithick] (-4,0) -- (4,0) node[right] {$x_1$};
\draw[ very thick] (0,0) circle (1.5);
\path (-3,0) [very thin,fill=black]  circle(1 pt) node at (-2.7, -0.2) { \fontsize{8}{8}\selectfont  {$ x=(x_1,0')$}};
\path (0,0) [very thin,fill=black]  circle(1 pt) node at (0, -0.2) { \fontsize{8}{8}\selectfont  {$ (1,0')$}};
\path (-1.5,0) [very thin,fill=black]  circle(1.7 pt) node at (-1.59, -0.2) { \fontsize{8}{8}\selectfont  {$ \textbf{0}$}};
\path  [very thin,fill=black]  circle(1 pt) node at (0, 1.2) {    {$ B$}};
\path  [blue,very thin,fill=black]  circle(1 pt) node at (1.8, 1.2) {    {$ D$}};
\draw[blue, very thick] (-1.5,-1.5) rectangle (1.5,1.5);
\node [below=1cm, align=flush center,text width=8cm] at (0,-0.6)
{$ Fig.2.$    \fontsize{10}{10}\selectfont  {Domain  D.} };
\end{tikzpicture}
\end{center}
We derive
\begin{eqnarray}\label{eq2-15}
 \int_{B_1(0)}\frac{1}{|x-y|^{n+s}}dy
&=&\int_{B_1(1, 0')}\frac{1}{|x-(y_1-1, y')|^{n+s}}dy\nonumber\\
&\leq& \int_{D}\frac{1}{|x-(y_1-1, y')|^{n+s}}dy\nonumber\\
&=& C \int_0^2 \int_0^1 \frac{r^{n-2}}{(r^2+(y_1-x_1-1)^2)^{\frac{n+s}{2}}}drdy_1\nonumber\\
&=& C \int_0^2 \frac{1}{|y_1-x_1-1|^{s+1}}dy_1\nonumber\\
& \sim & \frac{1}{|x_1|^s}.
\end{eqnarray}
Combining \eqref{eq2-14} with \eqref{eq2-15}, we derive
\eqref{eq2-16}.
\medskip

\emph{Step 2.}
In this step, we prove \eqref{eq2-mthm}.

For simplicity, we denote
$v(z)= \rho (x-z).$ Therefore, we need to show that
\begin{eqnarray*}
\int_{\mathbb{R}^n} (-\Delta)^s \phi(x) v(x)dx =\int_{\mathbb{R}^n} \phi(x)  (-\Delta)^s v(x)dx.
\end{eqnarray*}
According to the definition of the fractional Laplacian, the above equality is equivalent to
\begin{eqnarray}\label{eq2-17}
\int_{\mathbb{R}^n}  v(x)\lim_{\varepsilon \to 0}\int_{|y-x|\geq \varepsilon} \frac{\phi(x)-\phi(y)}{|x-y|^{n+2s}}dydx =\int_{\mathbb{R}^n} \phi(x)  \lim_{\varepsilon \to 0}\int_{|y-x|\geq \varepsilon}\frac{v(x)-v(y)}{|x-y|^{n+2s}}dydx.
\end{eqnarray}
Therefore, we only need to show \eqref{eq2-17} in this step.

We consider  the integral on the left-hand side and on the right-hand side of \eqref{eq2-17} respectively, we show that
the limit as $\varepsilon \to 0$ can be interchanged with the integral over $y.$

We first consider the integral on the left-hand side of \eqref{eq2-17} in $B^c_1(0)$ and $B_1(0)$ respectively.
If $x \in B_1^c(0),$ by \eqref{eq2-16}, we have
$$
\left| \int_{|y-x|\geq \varepsilon} \frac{\phi(x)-\phi(y)}{|x-y|^{n+2s}}dy \right| \leq \frac{C}{(dist(x, \partial B_1(0)))^s}.
$$
By Lebesgue's dominated convergence theorem, we obtain
\begin{eqnarray}\label{eq2-19}
&&\int_{B_1^c(0)} v(x) \lim_{\varepsilon \to 0}\int_{|y-x|\geq \varepsilon} \frac{\phi(x)-\phi(y)}{|x-y|^{n+2s}}dydx\nonumber\\
&=&\lim_{\varepsilon \to 0} \int_{B_1^c(0)}v(x) \int_{|y-x|\geq \varepsilon} \frac{\phi(x)-\phi(y)}{|x-y|^{n+2s}}dy dx.
\end{eqnarray}

If $x \in B_{1-\delta}(0),$ for any fixed $0<\delta <1$, we have
$$
\left| \int_{|y-x|\geq \varepsilon} \frac{\phi(x)-\phi(y)}{|x-y|^{n+2s}}dy \right| \leq C \|\phi\|_{C^{1,1}(B_{1-\delta}(0))}.
$$
Similarly, by Lebesgue's dominated convergence theorem, we derive
\begin{eqnarray*}
\int_{B_{1-\delta}(0)} v(x)\lim_{\varepsilon \to 0}\int_{|y-x|\geq \varepsilon} \frac{\phi(x)-\phi(y)}{|x-y|^{n+2s}}dy dx
=\lim_{\varepsilon \to 0} \int_{B_{1-\delta}(0)}v(x) \int_{|y-x|\geq \varepsilon} \frac{\phi(x)-\phi(y)}{|x-y|^{n+2s}}dy dx.
\end{eqnarray*}
Let $\delta \to 0,$ we obtain
\begin{eqnarray}\label{eq2-20}
\int_{B_1(0)} v(x)\lim_{\varepsilon \to 0}\int_{|y-x|\geq \varepsilon} \frac{\phi(x)-\phi(y)}{|x-y|^{n+2s}}dy dx
=\lim_{\varepsilon \to 0} \int_{B_1(0)}v(x) \int_{|y-x|\geq \varepsilon} \frac{\phi(x)-\phi(y)}{|x-y|^{n+2s}}dy dx.
\end{eqnarray}

Combining \eqref{eq2-19} with \eqref{eq2-20}, we derive that the limit $\varepsilon \to 0$ can be interchanged with the integral over $x$  on the left-hand side of \eqref{eq2-17} , i.e.,
\begin{eqnarray}\label{eq2-21}
\int_{\mathbb{R}^n}v(x) \lim_{\varepsilon \to 0}\int_{|y-x|\geq \varepsilon} \frac{\phi(x)-\phi(y)}{|x-y|^{n+2s}}dy dx
=\lim_{\varepsilon \to 0} \int_{\mathbb{R}^n}v(x) \int_{|y-x|\geq \varepsilon} \frac{\phi(x)-\phi(y)}{|x-y|^{n+2s}}dydx.
\end{eqnarray}

Then we consider the integral on the right-hand side of \eqref{eq2-17}. By Lebesgue's dominated convergence theorem, it is easy to see that the limit as $\varepsilon \to 0$ and the integral over $x$ on the right-hand side of \eqref{eq2-17} can be interchanged, i.e.,
\begin{eqnarray}\label{eq2-18}
\int_{\mathbb{R}^n} \phi(x)  \lim_{\varepsilon \to 0}\int_{|y-x|\geq \varepsilon}\frac{v(x)-v(y)}{|x-y|^{n+2s}}dydx
=\lim_{\varepsilon \to 0} \int_{\mathbb{R}^n} \phi(x)  \int_{|y-x|\geq \varepsilon}\frac{v(x)-v(y)}{|x-y|^{n+2s}}dydx.
\end{eqnarray}

We aim at showing  \eqref{eq2-17}, by \eqref{eq2-21} and \eqref{eq2-18}, it suffices to show that
\begin{eqnarray*}
 \int_{\mathbb{R}^n}v(x)\int_ {|y-x|\geq\varepsilon}\frac{\phi(x)-\phi(y)}{\mid x-y\mid^{n+2\sigma}}dy dx-\int_{\mathbb{R}^n} \phi(x)\int_ {|y-x|\geq\varepsilon}\frac{v(x)-v(y)}{\mid x-y\mid^{n+2\sigma}}dydx
 =0.
 \end{eqnarray*}
 In fact,  we derive from Fubini's Theorem that
\begin{eqnarray*}
 &&\int_{\mathbb{R}^n}v(x)\int_ {|y-x|\geq\varepsilon}\frac{\phi(x)-\phi(y)}{\mid x-y\mid^{n+2\sigma}}dy dx-\int_{\mathbb{R}^n} \phi(x)\int_ {|y-x|\geq\varepsilon}\frac{v(x)-v(y)}{\mid x-y\mid^{n+2\sigma}}dydx\\
 &=&\int_{\mathbb{R}^n}\int_ {|y-x|\geq\varepsilon}\frac{\phi(x)v(x)-\phi(y)v(x)}{\mid x-y\mid^{n+2\sigma}}dydx\\
 &&-\int_{\mathbb{R}^n}\int_ {|y-x|\geq\varepsilon}\frac{\phi(x)v(x)-\phi(x)v(y)}{\mid x-y\mid^{n+2\sigma}}dydx\\
 &=&\int_{\mathbb{R}^n}\int_ {|y-x|\geq\varepsilon}\frac{-\phi(y)v(x)}{\mid x-y\mid^{n+2\sigma}}dydx-\int_{\mathbb{R}^n} \int_ {|y-x|\geq\varepsilon}\frac{-\phi(x)v(y)}{\mid x-y\mid^{n+2\sigma}}dydx\\
 &=&-\int_{\mathbb{R}^n}\int_ {|y-x|\geq\varepsilon}\frac{\phi(y)v(x)}{\mid x-y\mid^{n+2\sigma}}dydx+\int_{\mathbb{R}^n} \int_ {|y-x|\geq\varepsilon}\frac{\phi(y)v(x)}{\mid x-y\mid^{n+2\sigma}}dydx\\
 &=&0.
 \end{eqnarray*}
 Therefore, \eqref{eq2-17} holds.
\medskip

\emph{Step 3.}
In this step, we prove \eqref{inter}.

By the definition of the fractional Laplacian and the mollifier, one has
\begin{eqnarray*}
&&( - \Delta )^s \phi_1 (x)\\
&=&C_{n,s} PV \int_{\mathbb{R}^n} \frac{\int_{\mathbb{R}^n}\rho(x-z)\phi(z)dz-\int_{\mathbb{R}^n}\rho(y-z)\phi(z)dz}{|x-y|^{n+2s}}  dy\\
&=& \int_{\mathbb{R}^n} \phi(z) (-\Delta)_x^s \rho(x-z)dz\\
&=& \int_{\mathbb{R}^n} \phi(z) (-\Delta)_z^s \rho(x-z)dz\\
&=& \int_{B_1(0)}  (-\Delta)_z^s \phi(z) \rho(x-z)dz+\int_{B_1^c(0)}  (-\Delta)_z^s \phi(z) \rho(x-z)dz\\
&\leq & \int_{B_1(0)}  (-\Delta)_z^s \phi(z) \rho(x-z)dz\\
&= &\int_{B_1(0)}  \lambda_1 \phi(z) \rho(x-z)dz\\
&= &\lambda_1 \phi_1 (x).
\end{eqnarray*}
Hence, we derive \eqref{inter}.

This completes the proof of Lemma \ref{mollification}.

\begin{remark}\label{mrem-inter}
It can be seen from the proof of the above lemma that if $u \in C^{1,1}_{loc} (\mathbb{R}^n)\cap   {\cal L}_{2s } $  and $v\in C_0^\infty(\mathbb{R}^n),$ we have
$$
\int_{\mathbb{R}^n}(-\Delta)^s u(x)\, v(x)dx=\int_{\mathbb{R}^n} u(x)\, (-\Delta)^s v(x)dx.
$$
\end{remark}

\end{document}